\documentclass[a4paper,psamsfonts,reqno]{amsart}
\usepackage{amscd,amssymb,amsmath}
\usepackage{amscd,amssymb}
\usepackage{verbatim}
\usepackage{hyperref}

\usepackage[all]{xy}

\newtheorem{Lemma}{Lemma}[section]
\newtheorem{Th}[Lemma]{Theorem}
\newtheorem{Prop}[Lemma]{Proposition}

\newtheorem{Cor}[Lemma]{Corollary}
\newtheorem{Def}[Lemma]{Definition}

\newtheorem{Rem}[Lemma]{Remark}

\newenvironment{Proof}{{\sc Proof.}\ }{~\rule{1ex}{1ex}\vspace{0.5truecm}}

\numberwithin{equation}{section}

\newcommand{\No}{{\mathbb N}_0}

\begin{document}
	
	\title{The Mittag-Leffler condition descents via pure monomorphisms}
	
\author{Dolors Herbera}

\thanks{Partially supported by the project  PID2020-113047GB-I00/AEI/10.13039/501100011033 financed by the State  Research Agency (Spain)}

	\address{Departament de Matem\`atiques \\
	Centre de Recerca Matem\`atica\\
		Universitat Aut\`onoma de Barcelona \\ 08193 Be\-lla\-te\-rra
		(Barcelona), Spain\\ e-mail: dolors@mat.uab.cat}
	
	\date{\today}
	
	\begin{abstract}
		This notes aims to clarify the proof given by Raynaud and Gruson \cite{RG} that the Mittag-Leffler property descents via pure rings monomorphism of commutative rings. A consequence of that is that projectivity dencents via such ring homomorphisms (cf.\cite{perry}), a  revision of the proof also allows to prove that the property of being pure-projective also descents via pure monomorphisms between commutative rings \cite[Section 8]{HPW}.
	\end{abstract}
	
	\maketitle

	In their fundamental paper \cite{RG}, Raynaud and Gruson introduced the class of Mittag-Leffler modules. They proved  how useful such notion was, showing an important number of striking results. One of them was the descent of projectivity via pure ring monomorphisms \cite[Th\'eor\`eme~II.3.1.3]{RG} (or universally injective maps, as they are named in \cite{RG}) of commutative rings.

	It seems there has been some misunderstanding in  the literature because, as noted  by Gruson in the  paper \cite{gruson},  statement \cite[Proposition~II.2.5.2]{RG} is wrong.  The descent of projectivity via pure monomorphisms is stated in  \cite[Examples~II.3.1.4]{RG} that are presented as a consequence of the wrong statement, and no correction for that is given in \cite{gruson}. However, to conclude the descent of projectivity via pure monomorphisms only \cite[Proposition~II.2.5.1, Th\'eor\`eme~II.3.1.3]{RG} are needed and these results are perfectly correct in the original paper.

The descent of projectivity means that if $R\to T$ is a pure ring monomorphism of commutative rings, and $M$ is a flat $R$-module then $M_R$ is projective if and only if $M\otimes _R T_T$ is a projective $T$-module. This result was reproved in \cite{perry} in the case the ring homomorphism $R\to T$ is faithfully flat. In  \cite{anger}, if was reproved for the case of pure-monomorphisms.  In both papers, it was also observed that results of Brewer and Rutter \cite[Theorem~2]{brewer-rutter} allow to state the result in the following way:

Let $R\to T$ be a pure ring monomorphism of commutative rings, and let $M$ be an $R$-module. Then $M_R$ is projective if and only if $M\otimes _R T_T$ is a projective $T$-module.

In a recent preprint, Herbera, Prihoda and Wiegand have observed that suitable  modifications of the original arguments due to Raynaud and Gruson allow also to show that pure projectivity descents via pure monomorphisms of commutative rings \cite{HPW}. The proof by Raynaud and Gruson is based on \cite[Proposition~II.2.5.1]{RG} which shows that the Mittag-Leffler propety descents via pure monomorphism of commutative rings. This result  is reproduced in  Proposition~\ref{descentpure}. 
	  
To make clear that the result is \cite{HPW} is correct, we have written this short note. In Section~\ref{schar} we introduce the characterization	of Mittag-Leffler modules that allows to prove \cite[Proposition~II.2.5.1]{RG}. The proof of the latter result is then included in Proposition~\ref{descentML}. Finally, in the third section we include the detailed proof of the descent of pure projectivity and, as a consequence, the descent of projectivity.

We stress the fact that we are just reproducing arguments that are already in \cite{RG}.

\medskip

We thank Michal Hrbek and Tom\'a\v s Lyson\v ek for pointing out an error in a too simple proof of  Lemma~\ref{descentflat}, and for providing reference \cite{master} that extends the results here to a suitable context of central extensions of rings, see Remark~\ref{further}.
	 
\section{Characterizations of Mittag-Leffler modules}\label{schar}

\begin{Def} \emph{\cite{RG}}
	Let $M$ be a right module over a ring $R$. Then $M$ is a
	\emph{Mittag-Leffler module} if the canonical map
	\[\rho \colon  M\bigotimes_R \prod _{i\in I}Q_i\to \prod _{i\in
		I}(M\bigotimes_RQ_i)\] is injective for any family $\{Q_i\}_{i\in
		I}$ of left $R$-modules.
\end{Def}

The following characterization of Mittag-Leffler modules is also due to Raynaud and Gruson. It is also reproved in \cite{perry}.

\begin{Prop} \label{charML} The following are equivalent conditions for a right $R$-module $M$.
		\begin{enumerate}
		\item [(i)] $M$ is a Mittag-Leffler module.
		\item[(ii)] Let $M=\varinjlim F_\alpha$ where $(F_\alpha, u_{\beta \alpha}\colon F_\alpha \to F_\beta)_{\alpha\le \beta \in \Lambda}$ is a directed system of finitely presented modules. Then for any $\alpha \in \Lambda$ there exists $\beta \ge \alpha$ such that, for any left $R$-module $Q$,  $\mathrm{ker} \,( u_{\beta \alpha}\otimes _R Q)= \mathrm{ker} \, (u_{\gamma \alpha}\otimes _RQ)$ for any $\gamma \ge \beta \in \Lambda$. 
		\item[(iii)] There exists a directed system of finitely presented modules $(F_\alpha, u_{\beta \alpha}\colon F_\alpha \to F_\beta)_{\alpha\le \beta \in \Lambda}$  such that $M=\varinjlim F_\alpha$, and satisfying that for any $\alpha \in \Lambda$ there exists $\beta \ge \alpha$ such that, for any left $R$-module $Q$,  $\mathrm{ker} \,( u_{\beta \alpha}\otimes _R Q)= \mathrm{ker} \,( u_{\gamma \alpha}\otimes _RQ)$ for any $\gamma \ge \beta \in \Lambda$.  
		
	\end{enumerate}
\end{Prop}

The key to prove Proposition~\ref{charML} is the following Lemma.

\begin{Lemma} \label{trick} Let $u\colon M\to N$ and $v\colon M\to M'$ be homomorphisms of right $R$-modules. Consider the push-out diagram
	
\[\begin{CD}
	M @>u>> N\\  @VvVV  @VVwV \\
	 M' @>u' >> B \end{CD}\]	
	
Then, there are exact sequences 
\[0\to \mathrm{ker}\, u \cap \mathrm{ker}\, v\hookrightarrow \mathrm{ker}\, u \stackrel{v}{\to} \mathrm{ker}\, u' \to 0  \]

\[0\to \mathrm{ker}\, u \cap \mathrm{ker}\, v\hookrightarrow \mathrm{ker}\, v \stackrel{u}{\to} \mathrm{ker}\, w \to 0  \]

In particular, $\mathrm{ker}\, u =\mathrm{ker}\, v$ if and only if $u'$ and $w$ are monomorphisms. 
\end{Lemma}

\begin{Proof} The proof is easily done using the push-out property combined with element-chasing.
	\end{Proof}

Using Lemma~\ref{trick}, the characterization of Proposition~\ref{charML} can be rewritten in the following fancy way:

\begin{Prop} \label{charML2} The following are equivalent conditions for a right $R$-module $M$.
	\begin{enumerate}
		\item [(i)] $M$ is a Mittag-Leffler module.
		\item[(ii)] Let $M=\varinjlim F_\alpha$ where $(F_\alpha, u_{\beta \alpha}\colon F_\alpha \to F_\beta)_{\alpha\le \beta \in \Lambda}$ is a directed system of finitely presented modules. Then for any $\alpha \in \Lambda$ there exists $\beta \ge \alpha$ such that,  for any $\gamma \ge \beta \in \Lambda$, the homomorphisms $w_{\gamma \alpha}$ and $u'_{\beta \alpha}$ in the push-out diagram
		
		\[\begin{CD}
			F_\alpha  @>u_{\beta \alpha}>> F_\beta \\  @Vu_{\gamma \alpha}VV  @VVw_{\gamma \alpha}V \\
			F_\gamma @>u'_{\beta \alpha} >> N_{\gamma \beta} \end{CD}\]	
	
		are pure monomorphisms. 
		
		\item[(iii)] There exists a directed system of finitely presented modules $(F_\alpha, u_{\beta \alpha}\colon F_\alpha \to F_\beta)_{\alpha\le \beta \in \Lambda}$  such that $M=\varinjlim F_\alpha$, and satisfying that for any $\alpha \in \Lambda$ there exists $\beta \ge \alpha$ such that,  for any $\gamma \ge \beta \in \Lambda$ the homomorphisms $w_{\gamma \alpha}$ and $u'_{\beta \alpha}$ in the push-out diagram
		
		\[\begin{CD}
			F_\alpha  @>u_{\beta \alpha}>> F_\beta \\  @Vu_{\gamma \alpha}VV  @VVw_{\gamma \alpha}V \\
			F_\gamma @>u'_{\beta \alpha} >> N_{\gamma \beta} \end{CD}\]	
		
	are pure monomorphisms.   
		\end{enumerate}
\end{Prop}

\section{Descent of the Mittag-Leffler condition via pure monomorphisms}

\begin{Prop} \label{descentpure}
	Let $\varphi \colon R\to T$ be a pure ring monomorphism of commutative rings. Let $M$, $N$ be $R$-modules  and let $u\colon M\to N$ an $R$-module homomorphism. Then $u$ is a pure monomorphism of $R$-modules if and only if $u\otimes _RT\colon M\otimes _RT\to N\otimes _RT$ is a pure monomorphism of $T$-modules. 
\end{Prop}

\begin{Proof}
	It is clear that if $u$ is a pure monomorphism then so is $u\otimes _RT$. 
	
	To prove the converse, assume that $u\otimes _RT$ is a pure monomorphism of $T$-modules. We claim it is also a pure monomorphism of $R$-modules.  Indeed, for any $R$-module $Q$ there is a commutative diagram
	
	\[\begin{CD}
		M\otimes _RT\otimes_R Q @>u\otimes _RT\otimes _RQ>> N\otimes _RT\otimes _RQ \\  @V\cong VV  @VV\cong 
		V \\
		M\otimes _RT\otimes _T(T\otimes_R Q) @>u\otimes _RT\otimes _T(T\otimes _RQ)>> N\otimes _RT\otimes _T(T\otimes _RQ)  \end{CD}\]
	in which the lower row is a monomorphism, then so is the upper row.
	
	Note also that for any $R$-module $X$, the embedding $X\otimes _R\varphi \colon X\to X\otimes _RT$ is a pure monomorphism of $R$-modules. 
	
	Finally,  the commutativity of the diagran
		\[\begin{CD}
		M @>M\otimes _R\varphi >> M\otimes _RT \\  @Vu VV  @VVu\otimes _RT
		V \\
		N @>N\otimes _R\varphi >> N\otimes _RT  \end{CD}\]
	in which $ (u\otimes _RT)\circ (M\otimes _R\varphi)$ is a pure monomorphism of $R$-modules, implies that $u$ is also a pure monomorphism of $R$-modules, as we wanted to prove.	
\end{Proof}

\begin{Lemma} \label{descentflat} \emph{(\cite{brewer-rutter}, \cite[Lemma~6]{anger})}
	Let $\varphi \colon R\to T$ be a pure ring monomorphism of commutative rings, and let $M_R$ be an $R$ module.  Then $M_R$ is flat if and only if the $T$-module $M\otimes _RT$ is flat.
\end{Lemma}

\begin{Prop} \label{descentML} \emph{(\cite[Proposition~II.2.5.1]{RG})}
	Let $\varphi \colon R\to T$ be a pure ring monomorphism of commutative rings, and let $M$ be an $R$-module. Then $M$ is a Mittag-Leffler $R$-module if and only if $M\otimes _RT$ is a Mittag-Leffler $T$-module.
\end{Prop}

\begin{Proof}
	Assume that $M_R$ is a Mittag-Leffler $R$-module. For any family of $T$-modules $\{Q_i\}_{i\in I}$ the canonical map $T\otimes _T\prod _{i\in I}Q_i\to \prod _{i\in I} T\otimes _TQ_i$ is an isomorphism. Hence,  the composition of maps
	$$(M\otimes _RT)\otimes_T \prod _{i\in I}Q_i\to M\otimes _R\left(  \prod _{i\in I} T\otimes _TQ_i \right) \to  \prod _{i\in I} M\otimes _R T\otimes _TQ_i$$
	is injective and, therefore, $M\otimes _RT$ is Mittag-Leffler as $T$-module.
	
	Now assume that $M\otimes _RT$ is a Mittag-Leffler $T$-module. Let $(F_\alpha, u_{\beta \alpha}\colon F_\alpha \to F_\beta)_{\alpha\le \beta \in \Lambda}$ is a directed system of finitely presented $R$-modules  such that $M=\varinjlim F_\alpha$. Then $(F_\alpha\otimes _RT, u_{\beta \alpha}\otimes _RT\colon F_\alpha \otimes _RT\to F_\beta \otimes _RT)_{\alpha\le \beta \in \Lambda}$ be a directed system of finitely presented $T$-modules  such that $M\otimes _RT=\varinjlim F_\alpha\otimes _RT$.
	
	By Proposition~\ref{charML}, for any $\alpha \in \Lambda$ there exists $\beta \ge \alpha$ such that, for any left $T$-module $Q$,  $\mathrm{ker} \, (u_{\beta \alpha}\otimes _RT \otimes _T Q)= \mathrm{ker} \, (u_{\gamma \alpha}\otimes _RT\otimes _TQ)$ for any $\gamma \ge \beta \in \Lambda$. In view of Proposition~\ref{charML2}, and since tensor products preserves push-out diagrams this is equivalent to say that for any $\alpha$ there exists $\beta \ge \alpha$ such that,  for any $\gamma$, in the push out diagram
		\[\begin{CD}
		F_\alpha  @>u_{\beta \alpha}>> F_\beta \\  @Vu_{\gamma \alpha}VV  @VVw_{\gamma \alpha}V \\
		F_\gamma @>u'_{\beta \alpha} >> N_{\gamma \beta} \end{CD}\]	
	$w_{\gamma \alpha}\otimes _RT$ and $u'_{\beta \alpha}\otimes _RT$ are pure monomorphisms of $T$-modules. By Proposition~\ref{descentpure}, we deduce that $w_{\gamma \alpha}$ and $u'_{\beta \alpha}$ are pure monomorphisms of $R$-modules. By Proposition~\ref{charML2}, we deduce that $M$ is a Mittag-Leffler $R$-module.	
	\end{Proof}

\section{Descent of pure projectivity}

A module $M$ is said to be pure projective if the functor $\mathrm{Hom}_R(M, -)$ is exact with pure short exact sequences. Equivalently, $M$ is pure projective if it is a direct summand of a direct sum of finitely presented modules.

Pure-projective modules always decompose into a direct sum of countably presented pure-projective submodules. 

\medskip

Now we include the proof that pure projectivity descents via pure monomorphisms. We reproduce this result from \cite[\S 8]{HPW}.

First we recall the results that relate pure-projective modules and Mittag-Leffler modules.

\begin{Rem}\label{rem:walrus}  The map $\rho$ in the definition of Mittag-Leffler module is obviously bijective if $M$ is a finitely
	generated free module. An easy diagram chase shows that it is also bijective if 
	$M$ is finitely presented.  Thus finitely presented modules are Mittag-Leffler.  Since the class of 
	Mittag-Leffler modules is closed under direct summands and arbitrary direct sums,  all pure-projective modules are Mittag-Leffler modules.
\end{Rem}

\begin{Lemma} \label{ML2} \cite[Th\'eor\`eme~2.2.1 p. 73]{RG} 
	Any countably generated submodule $X$  of a Mittag-Leffler module $Y$ is contained in a pure-projective countably generated  pure submodule $Y'$ of $Y$. 
\end{Lemma}

The following Proposition is a variation of \cite[Th\'eor\`eme~3.1.3 p. 78] {RG} adapted to the pure-projective situation. 

\begin{Prop} \emph{\cite{HPW}} \label{descentpureproj} Let $R\subseteq T$ be a pure extension of commutative rings. Let $M$ be an $R$-module. Then $M_R$ is pure projective  if and only if $M\otimes _RT$ is pure projective as a $T$-module.
\end{Prop}

\begin{Proof}  If $M_R$ is pure projective then, clearly,   $M\otimes _RT$ is pure projective as a $T$-module. 
	For the converse, assume that $M\otimes _RT$ is a pure projective $T$-module.  By Remark~\ref{rem:walrus} $M\otimes_RT$ is a Mittag-Leffler $T$-module, and then Proposition~\ref{descentML} implies that $M$ is a Mittag-Leffler $R$-module.  We need to prove that, in addition, $M_R$ is pure-projective.

	Write  $M\otimes _RT=\oplus _{i\in I}Q_i$, where each $Q_i$ is a countably generated $T$-module.  Let $\mathcal F = \{Q_i\}_{i\in I}$.
	A submodule $X$ of $M$ is said to be \emph{adapted} (to $\mathcal{F}$) if it is pure and the canonical  image of $X\otimes _RT$ in $M\otimes _RT$ is a direct sum of modules in $\mathcal{F}$. If $X$ is an adapted submodule of $M$, then the sequence
	\[
	0\to X\otimes _RT\to M\otimes _RT\to (M/X)\otimes _RT \to 0
	\]
	is split exact. Therefore 
	$X\otimes _RT$  and $(M/X)\otimes _RT$, being isomorphic to direct summands of $M\otimes_RT$, are
	pure-projective as $T$-modules. 
	
	\smallskip
	
	\noindent\textbf{Step 1.} \emph{Every countably generated pure submodule of $M$ is contained in a countably generated adapted submodule of $M$. }
	
	Let $X$ be a countably generated pure submodule of $M$. As $M$ is Mittag-Leffler and the modules $Q_i$ are countably generated, we can construct
	a sequence $(X_n, I_n)_{n\in \No}$  such that
	\begin{itemize}
		\item[(1)] $X_0=X$ and, for every $n\ge 0$, $X_n$ is a countably generated pure submodule of $M$ and $X_n\subseteq X_{n+1}$;
		\item[(2)]  for any $n\ge 0$, $I_n$ is a countable subset of $I$ and it consists of the elements $i\in I$ such that the canonical projection of $X_n\otimes _RT$  in $Q_i$ is different from zero;
		\item[(3)]  for any $n\ge 0$, the image of $X_{n+1}\otimes _RT$ contains $\oplus _{i\in I_n}Q_i$.
	\end{itemize} 
	To be more specific, 
	suppose $X_n,I_n$ were defined. Since each $Q_i$ is countably generated, there exists 
	a countable set $G_n \subseteq M$ such that the canonical image of $G_nR \otimes_R T$ in 
	$M \otimes_R T$ contains $\oplus_{i \in I_n} Q_i$. By Lemma~\ref{ML2}, there exists 
	a countably generated $X_{n+1}$ which is a pure submodule of $M$ containing $X_n + G_nR$.
	Then $I_{n+1}$ is chosen as described in (2).
	
	Set $Y=\bigcup _{n\in \No}X_n$. By construction, $Y$ is an adapted submodule of $M$.
	
	\smallskip
	
	\noindent\textbf{Step 2.} \emph{ Let $X$ be an arbitrary adapted submodule of $M$ such that $X\neq M$. 
		Then there exists an adapted 
		submodule $X'$ of $M$
		such that $X\subset X'$ and $X'/X$ is a countably generated adapted 
		submodule of $M/X$. Hence $X'/X$ is pure-projective; therefore the pure exact sequence
		\[
		0\to X\to X'\to X'/X\to 0
		\]
		splits.}
	
	By definition, if $X$ is an adapted submodule of $M$ then $(M/X)\otimes _RT\cong \oplus _{i\in I'}Q_i$  for a certain $I'\subseteq I$. Hence,  it makes sense to talk about adapted submodules of $M/X$ with respect to the decomposition induced by $\mathcal{F}'=\{Q_i\}_{i\in I'}$. By Step 1, there exists a submodule $X'$ of $M$ containing $X$ and
	such that $X'/X$ is a countably generated adapted submodule of $M/X$. Therefore  $X'$ is also an adapted submodule of $M$. Since $X$ is an adapted submodule of $X'$ the rest of the statement is clear.
	
	\smallskip
	
	Finally, combining the first and the second steps, we deduce that there exist an ordinal $\kappa $ and a continuous chain $\{X_\alpha\}_{\alpha <\kappa}$ of adapted submodules of $M$ such that
	\begin{itemize}
		\item[(i)] $X_0=0$, and
		\item[(ii)] for any $\alpha +1 <\kappa$, $X_{\alpha +1}/X_{\alpha}$ is pure projective and a direct summand of $X_{\alpha +1}$. 
	\end{itemize}
	In this situation $M \cong \oplus (X_{\alpha +1}/X_{\alpha})$ (cf. \cite[Lemme~3.1.2, p. 81]{RG}). Therefore $M$ is pure projective.
\end{Proof}

\begin{Cor} \label{descentprojective}
	Let $R\to T$ be a pure ring monomorphism of commutative rings, and let $M$ be an $R$-module. Then $M_R$ is projective if and only if $M\otimes _R T$ is a projective $T$-module.
\end{Cor}

\begin{Proof}
	Since projective modules are exactly the pure-projective modules that, in addition, are flat, the statement follows from Proposition~\ref{descentpureproj} and Lemma~\ref{descentflat}.
\end{Proof}

\begin{Rem} \label{further}  All the arguments presented, make heavy use of the commutativity hypothesis. It is worth to remark that Tom\'a\v s Lyson\v ek in his Master Thesis \cite{master} has made a detailed analysis of the descent properties via homomorphisms of general associative rings.
	
For example in \cite[Corrollary~7.6]{master} he gives a new approach to the descent of the property of being flat via pure ring monomorphisms that we consider in Lemma~\ref{descentflat}.	He also proves the following interesting result, that extends the ones proved in this note and  allows some     non-commutativity:
	
	\begin{Th} \cite[Theorem~8.5]{master}
	Let $R$ be a commutative ring, $\varphi \colon R\to S$ be a central pure ring homomorphism, let $A$ be a central $R$-algebra.  Then projectivity and pure-projectivity of both left and right modules descends through the $R$-algebra homomorphism $\varphi _A\colon A\to A\otimes _RS$ defined by  $\varphi _A (a)=a\otimes \varphi(1)$.
	\end{Th}	
	
\end{Rem}

\end{document}